\documentclass[a4paper, 12pt]{article}

\usepackage{graphicx}
\usepackage{epsfig}

\newtheorem{theorem}{Theorem}
\newtheorem{lemma}{Lemma}
\newtheorem{remark}{Remark}
\newtheorem{problem}{Problem}
\newtheorem{proposition}{Proposition}
\newtheorem{conjecture}{Conjecture}
\newtheorem{corollary}{Corollary}

\title{Semidefinite Programming and Reachable Sets of Dissipative
Bilinear Control Systems}

\author{\\Dionisis Stefanatos\quad and \quad Navin Khaneja\thanks{
D. Stefanatos and N. Khaneja are with the Division of Engineering
and Applied Sciences, Harvard University, Cambridge, MA 02138 USA
(e-mail: stefanat@fas.harvard.edu; navin@eecs.harvard.edu).}}

\date{}

\begin{document}

\maketitle  \pagestyle{plain}

\begin{abstract}
In this manuscript, we investigate optimal control problems which
arise in connection with manipulation of dissipative quantum
dynamics. These problems motivate the study of a class of
dissipative bilinear control systems. For these systems it is
shown that the optimal solution and the reachable set can be found
by solving a semidefinite program. In practice, solutions to these
problems generate optimal methods for control of quantum
mechanical phenomena in presence of dissipation. In the area of
coherent spectroscopy, this translates into the maximum signal to
noise ratio that can be obtained in a spectroscopy experiment.
\end{abstract}

\section{Introduction-Statement of the Problem}

Consider the following optimal control problem. Given the
dynamical system below
\begin{equation}
\label{physicalsystem} \left[\begin{array}{cccc}
\dot{x}_1\\\dot{x}_2\\\dot{y}_1\\\dot{y}_2\end{array}\right]=
\left[\begin{array}{cccc} 0 & 0 & -v_1 & 0
\\0 & 0 & 0 & -v_2\\ v_1 & 0 &-k & -J\\ 0 & v_2 & J & -k \end{array}\right]
\left[\begin{array}{cccc}x_1\\x_2\\y_1\\y_2\end{array}\right]
\end{equation}
and starting from the initial state $e_1=(1,0,0,0)^T$, what is the
maximum achievable value of $x_2$ and what are the optimal
controls $v_1(t) \in \Re$ and $v_2(t) \in \Re$ that achieve this
value? Problems like this are associated with optimal manipulation
of quantum mechanical phenomena under dissipation. Specifically,
the optimization problem stated above comes from Nuclear Magnetic
Resonance (NMR) spectroscopy and is related to optimal control of
two coupled spins in presence of transverse relaxation
\cite{Khaneja03}. The state variables $x_i, y_i$ represent
averages of various quantum mechanical spin operators. The
available controls $v_1(t)$ and $v_2(t)$ correspond to the
components of the magnetic field in the NMR experimental setup.
Parameter $k>0$ expresses the transverse relaxation rate while $J$
is the coupling constant between the spins.

\begin{figure}[h]
\label{Fig1} \centering
\includegraphics[scale=0.4]{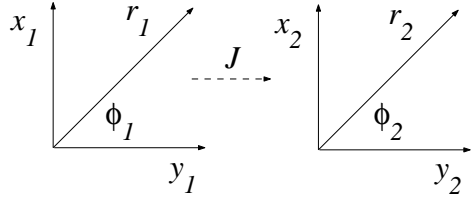}
\caption{\label{r1r2}Schematic representation of the evolution of
system (\ref{physicalsystem}). Control $v_1$ rotates $x_1$ to
$y_1$. Under the $J$ coupling $y_1$ evolves to $y_2$, while both
dissipate because of the relaxation term $k$. Control $v_2$
rotates $y_2$ to $x_2$. The new state variables $r_{1}, r_{2}$,
defined by (\ref{r}), are also shown. The corresponding new
control parameters are $u_1=\cos\phi_1, u_2=\cos\phi_2$.}
\end{figure}

Observe that if $v_1$ and $v_2$ are set to 0 then the initial
state $e_1$ doesn't evolve at all and there is no build up of
$x_2$. However, by turning on $v_1$, it is possible to rotate
$x_1$ to $y_1$, see Fig. \ref{r1r2}. This evolves to $y_2$ under
the skew symmetric matrix
\begin{displaymath}
\left[\begin{array}{cc}0 & -J
\\J & 0\end{array}\right]\;,
\end{displaymath}
while both $y_1$ and $y_2$ dissipate under the term
\begin{displaymath}
\left[\begin{array}{cc}-k & 0
\\0 & -k\end{array}\right]\;.
\end{displaymath}
The state $y_2$ can then be rotated to $x_2$ by switching on the
control $v_2$. We want to find the optimal $v_1$ and $v_2$ that
maximize the value of $x_2$. It is intuitively clear that no
matter how large we make $v_1(t), v_2(t)$, the transfer
$x_1\rightarrow x_2$ cannot be done without any loss, since the
intermediate transfer $y_1\rightarrow y_2$ is entirely due to
internal dynamics over which there is no control, thus there is an
unavoidable dissipation because of $k>0$.

Define
\begin{equation}
\label{r}r_i=\sqrt{x_i^2+y_i^2}\;.
\end{equation}
Using (\ref{physicalsystem}), evolution equations for $r_1, r_2$
can be found. We get the system
\begin{displaymath}
\left[\begin{array}{cc} \dot{r}_{1}\\\dot{r}_{2}\end{array}\right]
= \left[\begin{array}{cc}-k\cos^2\phi_1 & -J
\cos{\phi_1}\cos{\phi_2}
\\J\cos{\phi_1}\cos{\phi_2} &
-k\cos^2\phi_2\end{array}\right]\left[\begin{array}{cc}
r_{1}\\r_{2}\end{array}\right],
\end{displaymath}
where $\cos{\phi_1}=y_1/r_1, \cos{\phi_2}=y_2/r_2$, see Fig.
\ref{r1r2}. Using the control $v_1$ which rotates $x_1$ to $y_1$,
we can control the angle $\phi_1$. Analogously, using $v_2$ which
rotates $y_2$ to $x_2$, we can control the angle $\phi_2$.
Denoting $u_1=\cos{\phi_1}$, $u_2=\cos{\phi_2}$ and dilating time
by a factor of $J$, the above system can be rewritten as
\begin{equation}
\label{system2} \left[\begin{array}{c}
\dot{r}_{1}\\\dot{r}_{2}\end{array}\right] =
\left[\begin{array}{cc}-\xi u_{1}^{2} & -u_{1}u_{2}
\\ u_{2}u_{1} & -\xi u_{2}^{2}\end{array}\right]\left[\begin{array}{c}
r_{1}\\r_{2}\end{array}\right]\;.
\end{equation}
Here $u_1$ and $u_2$ are control parameters which take their
values in the interval $[-1, 1]$ and $\xi=k/J$. The initial
problem of maximum transfer from $x_1$ to $x_2$ has been
transformed to the following equivalent question:

Given the dynamical system (\ref{system2}) and the initial state
$(r_1(0),r_2(0))= (1,0)$, find the optimal control
$(u_1(t),u_2(t))$, $|u_1|,|u_2|\leq 1$, such that $r_2$ is
maximized.

Note that once $r_2$ is maximized, the control $v_2$ can be used
to transfer it to $x_2$ with no loss, so the above question is
indeed equivalent to the original problem.

Motivated by this example, which originates from a real physical
system, let us consider the following $n$-dimensional
generalization of system (\ref{physicalsystem}):
\begin{equation}
\label{systemn} \left[\begin{array}{c}
\dot{x}\\\dot{y}\end{array}\right] = \left[\begin{array}{cc} 0 &
-V
\\ V & A\end{array}\right]\left[\begin{array}{c}
x\\y\end{array}\right]\;,
\end{equation}
where $x=(x_1,x_2,\ldots,x_n)^T$, $y=(y_1,y_2,\ldots,y_n)^T$,
$V=\mbox{diag}(v_1,v_2,\ldots,v_n)$ and $A=\{a_{ij}\}$ is such
that its symmetric part $A+A^T$ is negative definite. This
condition insures that the norm of the vector $(x,y)$ can only
decrease. This models the physics in open quantum systems, where
dissipation can only reduce coherence in the system dynamics.
Furthermore,  $A$ is such that any two states $y_i$ and $y_j$ are
coupled by its off-diagonal elements, not necessarily directly (we
say $A$ is \emph{irreducible}).
\begin{problem}
\label{bilinear}Given the dynamical system (\ref{systemn}) and the
starting state $(x(0),y(0))$, find the optimal control
$(v_1(t),v_2(t),\ldots,v_n(t))$ which maximizes $x_n$.
\end{problem}

It is shown in the following section that the negative
definiteness condition on $A$ is a sufficient condition for the
existence of an optimal solution.

If we define $r_i=\sqrt{x_i^2+y_i^2}$ and work as in the
$2$-dimensional case, we find that $r_i$ satisfies the equation
\begin{equation}
\label{generalsystem}\frac{dr_i}{dt}=\sum_{j=1}^na_{ij}u_iu_jr_j\;,
\end{equation}
where $u_i=y_i/r_i$. Problem \ref{bilinear} has been transformed
to the following.
\begin{problem}
\label{problem}Given the dynamical system defined by
(\ref{generalsystem}) for $i=1,2,\ldots,n$, with $A=\{a_{ij}\}$
irreducible and such that $A+A^T$ negative definite, and the
starting state $(r_1(0),r_2(0),\ldots,r_n(0))$, with $r_i(0)\geq
0$, find the optimal control $(u_1(t),u_2(t),\ldots,u_n(t))$,
$|u_i|\leq 1$, which maximizes $r_n$, while it preserves
$r_i(t)\geq 0$.
\end{problem}

Observe that if $T_1<T_2$, then the maximum achievable value in
time $T_1$ cannot exceed the corresponding value in time $T_2$,
since by putting $u_i=0$ the evolution in the interval $(T_1,T_2]$
can be stopped. Therefore, Problem \ref{problem} is considered as
an \textit{infinite horizon problem}.

Multiplying the $i^{th}$ equation of system (\ref{generalsystem})
with $2r_i$, we get
\begin{equation}
\frac{d}{dt}\left(r_i^2\right)=\sum_{j=1}^n2a_{ij}u_iu_jr_ir_j
\end{equation}
and from this
\begin{equation}
\label{ri2}
\frac{d}{dt}\left(r_i^2\right)=U^2\sum_{j=1}^n2a_{ij}\frac{u_ir_i}{U}\frac{u_jr_j}{U}\;,
\end{equation}
where
\begin{equation}
U = \sqrt{\sum_{i=1}^n (u_i r_i)^2}\;.
\end{equation}
By setting
\begin{equation}
p_i=r_i^2\,,\quad\quad m_i=\frac{u_ir_i}{U}\;,
\end{equation}
and rescaling time according to $dt'=U^2dt$, equation (\ref{ri2})
becomes
\begin{equation}
\label{psystem} \frac{dp_i}{dt'}=\sum_{j=1}^n2a_{ij}m_im_j\;.
\end{equation}
The initial optimal control problem has been transformed to the
following one.
\begin{problem}\label{mproblem}
Given the dynamical system defined by (\ref{psystem}) for
$i=1,2,\ldots,n$ and the starting point
$p(0)=(p_1(0),p_2(0),\ldots,p_n(0))^T$, $p_i(0)\geq$0, find the
unit vector $m(t')=(m_1(t'),m_2(t'),\ldots,m_n(t'))^T$ that
maximizes $p_n$, while it preserves $p_i(t')\geq 0$. Matrix
$A=\{a_{ij}\}$ is irreducible and such that $A+A^T$ is negative
definite.
\end{problem}

Note that, although Problem \ref{problem} is an infinite horizon
problem, Problem \ref{mproblem} defined above may achieve its
maximum for a finite final time $T_f$. There is no inconsistency
here, since the times for the two systems are related through
$dt'=U^2dt$, so $T_f=\int_0^{T_f}dt'=\int_0^\infty U^2dt$. If
$U(t)\rightarrow 0$ sufficiently fast as $t\rightarrow\infty$,
then $T_f$ is finite. As we will see, this is indeed the case.

In the following, we study problems \ref{problem} and
\ref{mproblem} in detail. Having found an optimal solution for the
latter, we can easily find a corresponding optimal control law for
the former. The structure of the paper is as follows. In section
\ref{SolStrat}, it is shown that the solution of Problem
\ref{mproblem} can be reduced to the solution of a semidefinite
program and that the negative definiteness of the symmetric part
of $A$ is a sufficient condition for the existence of an optimal
solution. It is also shown how the semidefinite programming
formalism can be used for calculating reachable sets. In section
\ref{Ranks}, some useful existent results regarding the rank of
matrices that solve our semidefinite program are presented. These
results are used in section \ref{Examples}, where some specific
examples are examined. The examples include system (\ref{system2})
and another system which again arises from an optimal control
problem of spin dynamics in NMR spectroscopy.

\section{\label{SolStrat}Reduction to a Semidefinite Program}

In the following, the inner product $\langle\cdot\,,\cdot\rangle$
in the space of symmetric $n\times n$ matrices $\mbox{Sym}_n$ is
defined in the usual way as the trace of the matrix product, i.e.
$\langle A, B\rangle=\mbox{tr}(AB)$ for $A,B\in\mbox{Sym}_n$. Note
also that $A\succeq 0$ denotes that matrix $A\in\mbox{Sym}_n$ is
positive semidefinite, $A\prec 0$ that is negative definite etc.
\begin{theorem}
\label{semidefinite}Let us define matrices $A_i\in\mbox{Sym}_n,
i=1,2,\ldots,n$, by the relation
\begin{displaymath}
A_i=\left[\begin{array}{ccccccc}
 & & &a_{i1}& & & \\
 &\textbf{O}_{i1}& &\vdots& &\textbf{O}_{i2}& \\
 & & &a_{i(i-1)}& & & \\
a_{i1}&\ldots&a_{i(i-1)}&2a_{ii}&a_{i(i+1)}&\ldots&a_{in} \\
 & & &a_{i(i+1)}& & & \\
 &\textbf{O}_{i3}& &\vdots& &\textbf{O}_{i4}& \\
 & & &a_{in}& & &
\end{array}
\right]\;,
\end{displaymath}
where $a_{ij}$ are the elements of matrix $A$ given in Problem
\ref{mproblem} and $\textbf{O}_{il}, l=1,2,3,4$, are zero matrices
with appropriate sizes. The solution of Problem \ref{mproblem} can
be reduced to the solution of the following semidefinite program:
\begin{displaymath}
\mbox{Find}\quad\mathcal{E}=\max_M\langle A_n,M\rangle
\end{displaymath}
\begin{displaymath}
\mbox{subject to}\quad\langle A_i,M\rangle=-p_i(0)\,,\quad
i=1,2,\ldots,n-1
\end{displaymath}
\begin{displaymath}
\mbox{and}\quad M\succeq 0\;.
\end{displaymath}
The maximum achievable value of $p_n$ is $p_n(0)+\mathcal{E}$.
\end{theorem}
\begin{proof}
Let $T_f$ be the time when $p_n$ achieves its maximum, i.e. the
final time. From equation (\ref{psystem}) it is
\begin{equation}
\label{pint}
p_i(T_f)=p_i(0)+\sum_{j=1}^n2a_{ij}\int_0^{T_f}m_i(t')m_j(t')dt'\;.
\end{equation}
Observe that if we define the positive semidefinite matrix $M$
through the relation
\begin{equation}
\label{M} M=\int_0^{T_f}m(t')m^T(t')dt'\;,
\end{equation}
then (\ref{pint}) becomes
\begin{equation}
\label{pin} p_i(T_f)=p_i(0)+\langle A_i,M\rangle\;.
\end{equation}

One other important observation is that the end point of the
optimal trajectory should lie on the line $(0,0,\ldots,0,p_n)$ in
$p$-space. Suppose that the end point has a component $p_k>0$ for
some $k\neq n$. If $p_k$ is directly coupled to $p_n$ then choose
$m=(0,0,\ldots,0,m_k,0,\ldots,0,m_n)^T$ such that
$m_n(a_{nk}m_k+a_{nn}m_n)>0$ and evolve the system until $p_k=0$.
Thereby we get a greater value of $p_n$. If $p_k$ is not directly
coupled to $p_n$, we can still transfer from $p_k$ to $p_n$ using
intermediate states (because matrix $A$ is irreducible). We
conclude that at the final time $T_f$ the end point of the optimal
trajectory should lie on the line $(0,0,\ldots,0,p_n)$. Thus, we
have to maximize $p_n(T_f)=p_n(0)+\langle A_n,M\rangle$ under the
conditions $p_i(T_f)=p_i(0)+\langle A_i,M\rangle=0$,
$i=1,2,\ldots,n-1$. Equivalently, we have to solve the following
semidefinite program: Find $\mathcal{E}=\max_M\langle
A_n,M\rangle$ subject to $\langle A_i,M\rangle=-p_i(0)$ for
$i=1,2,\ldots,n-1$ and $M\succeq 0$.

Having found an optimal $M$, we can always find an appropriate
unit vector $m(t')$ such that $M=\int_0^{T_f}m(t')m^T(t')dt'$ and
$p_i(t')\geq 0$. Since $M\succeq 0$, it can always be decomposed
in the form
\begin{displaymath}
M=\sum_{k=1}^r\lambda_k\textbf{m}_k\textbf{m}_k^T\;,
\end{displaymath}
where $\lambda_k$ are the positive eigenvalues of $M$,
$\textbf{m}_k$ are the corresponding (real) normalized
eigenvectors and $r$ is the rank of $M$. Now let $N$ be a positive
integer. Rewrite the above relation in the form
\begin{displaymath}
M=N\sum_{k=1}^r\Delta\lambda_k\textbf{m}_k\textbf{m}_k^T\;,
\end{displaymath}
where $\Delta\lambda_k=\lambda_k/N$, and define the times $t_k'$
through
\begin{displaymath}
t_0'=0\;,\quad
t_k'=\sum_{l=1}^k\Delta\lambda_l\quad\mbox{for}\quad
k=1,2,\ldots,r\;.
\end{displaymath}
Let us forget for a moment the restrictions $p_i(t')\geq 0$. If we
apply the control
\begin{displaymath}
m(t')=\textbf{m}_k\quad\mbox{for}\quad t_{k-1}'\leq t'<t_k'\,,
\quad k=1,2,\ldots,r
\end{displaymath}
and repeat for $N$ times, then on the one hand the requirement
\begin{displaymath}
\int_0^{T_f}m(t')m^T(t')dt'=M
\end{displaymath}
is satisfied and on the other hand the trajectory in $p$-space
approximates the line joining the initial point
$I(p_1(0),p_2(0),\ldots,p_n(0))$ to the final point
$F(0,0,\ldots,p_n(T_f))$, see Fig \ref{ptrajectory}(a). If $N$ is
large enough then the trajectory actually follows this line, see
Fig. \ref{ptrajectory}(b), thus the restrictions $p_i(t')\geq 0$
are satisfied. Note that $T_f=\sum_{k=1}^r\lambda_k=\mbox{tr}(M)$
is finite, if $\mbox{tr}(M)<+\infty$. In the special case where
$r=1$, it is $M=\lambda\textbf{m}\textbf{m}^T$ and thus
$m(t')=\textbf{m}$ for $t'\in [0,T_f]$, $T_f=\lambda$.

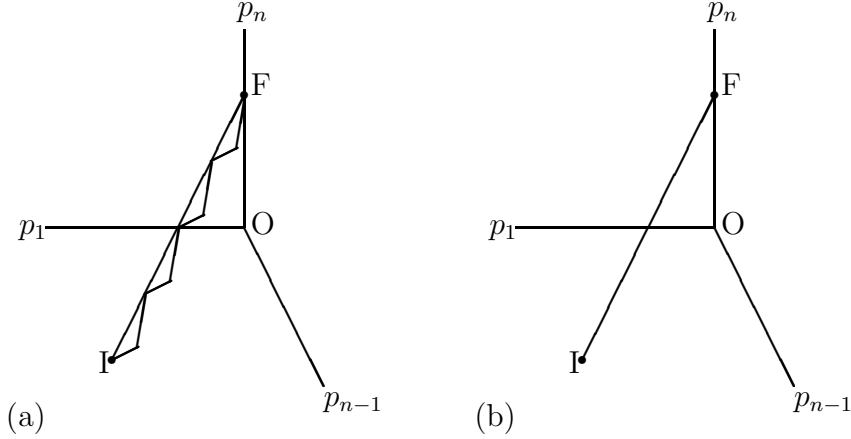
\begin{figure}
\setlength{\unitlength}{0.88cm} \centering
\begin{tabular}{cc}
(a)\begin{picture}(6,6) \thicklines \put(3,3){\line(0,1){3}}
\put(3,3){\line(-1,0){3}}\put(3,3){\line(1,-2){1.2}}\put(1,1){\line(1,2){2}}
\put(3.1,2.9){O}\put(2.9,6.2){$p_n$}\put(-0.4,2.9){$p_1$}\put(4.2,0.3){$p_{n-1}$}\put(0.8,0.8){I}\put(1,1){\circle*{0.1}}\put(3.1,5){F}\put(3,5){\circle*{0.1}}
\multiput(1,1)(0.5,1){4}{\line(2,1){0.4}}\multiput(1.38,1.1818)(0.5,1){4}{\line(1,6){0.135}}
\end{picture}
& (b)\begin{picture}(6,6) \thicklines \put(3,3){\line(0,1){3}}
\put(3,3){\line(-1,0){3}}\put(3,3){\line(1,-2){1.2}}\put(1,1){\line(1,2){2}}
\put(3.1,2.9){O}\put(2.9,6.2){$p_n$}\put(-0.4,2.9){$p_1$}\put(4.2,0.3){$p_{n-1}$}\put(0.8,0.8){I}\put(1,1){\circle*{0.1}}\put(3.1,5){F}\put(3,5){\circle*{0.1}}
\end{picture}
\end{tabular}
\caption{\label{ptrajectory}(a) Trajectory in $p$-space following
the control law presented in the text, for $r=2$ and $N=4$. It
approximates the straight line from the initial point $I$ to the
final point $F$. Note that the restrictions $p_i(t')\geq 0$ may
not be satisfied when $N$ is small (b) For large $N$ the
trajectory coincides with the line $IF$, so the restrictions
$p_i(t')\geq 0$ are satisfied.}
\end{figure}

The conclusion is that we just need to solve the semidefinite
program defined above. The maximum achievable value of $p_n$ is
$p_n(T_f)=p_n(0)+\mathcal{E}$.
\end{proof}

We show next how this control law can be applied to system
(\ref{generalsystem}) in Problem \ref{problem}. For $0\leq t'\leq
t_1'$, $m(t')=\textbf{m}_1=\mbox{constant}$. Since, additionally,
$\textbf{m}_1$ is a unit vector, we can assume without loss of
generality that its first component $m_1\neq 0$. Consider the
ratios
\begin{displaymath}
\frac{u_ir_i}{u_1r_1}=\frac{m_i(t')}{m_1(t')}=s_i\,,\quad\;i=1,2,\ldots,n.
\end{displaymath}
For $0\leq t'\leq t_1'$, $s_i$ are constant. Define
\begin{displaymath}
\mathcal{M}=\max_i\left(\left|\frac{s_ir_1}{r_i} \right|
\right)\,,\quad i=1,2,\ldots,n.
\end{displaymath}
The optimal policy can be realized as
\begin{displaymath}
u_1=\frac{1}{\mathcal{M}}
\end{displaymath}
and
\begin{displaymath}
u_i =\frac{s_ir_1}{r_i}u_1\;,
\end{displaymath}
where $i=2,3,\ldots,n$. With the above choice we insure that
$|u_i|\leq 1$. Using this feedback law we can evolve system
(\ref{generalsystem}) in time $t$ and calculate the function
$U(t)=\sum_{i=1}^n(u_ir_i)^2$. Then, we can find
$t'=\int_0^tU^2dt$. When $t'=t_1'$, we switch to
$m(t')=\textbf{m}_2$ and repeat the above procedure. If the rank
of $M$ is $r=1$ then the ratios $s_i$ keep the same value for all
times. Note that the maximum achievable value of $r_n$ is
\begin{equation}
r_n(\infty)=\sqrt{p_n(T_f)}=\sqrt{p_n(0)+\mathcal{E}}=\sqrt{r_n^2(0)+\mathcal{E}}\;.
\end{equation}

In the above discussion we implicitly assumed that an optimal
solution exists, and we used for $\mathcal{E}$ the
characterization ``maximum" instead of the more formal
``supremum". We show below that the negative definiteness of
$A+A^T$ is a sufficient condition for the existence of an optimal
solution. The following lemma is used.
\begin{lemma}
\label{BM}If $B\succ 0$ and $M\succeq 0$, $B,M\in\mbox{Sym}_n$,
then $\langle B,M\rangle\geq 0$.
\end{lemma}
\begin{proof}
Since $B\in\mbox{Sym}_n$ it can be diagonalized by an orthogonal
matrix $O$, $B=O\Delta O^T$, where
$\Delta=\mbox{diag}(\lambda_1,\lambda_2,\ldots,\lambda_n)$ and
$\lambda_i>0$ are the eigenvalues of the positive definite matrix
$B$. It is
\begin{displaymath}
\langle B,M\rangle=\mbox{tr}(BM)=\mbox{tr}(O\Delta
O^TM)=\mbox{tr}(\Delta O^TMO)=\mbox{tr}(\Delta
\tilde{M})=\sum_{i=1}^n\lambda_i\tilde{m}_{ii}\;,
\end{displaymath}
where $\tilde{M}=O^TMO$ and $\tilde{m}_{ii}$ are its diagonal
elements. But $\tilde{M}^T=\tilde{M}$ and
$x^T\tilde{M}x=(Ox)^TMOx\geq 0$ for every $x\in\Re^n$, since
$M\succeq 0$. So, it is also $\tilde{M}\succeq 0$ and thus
$\tilde{m}_{ii}\geq 0$. Since, additionally, $\lambda_i>0$, we
conclude that $\langle
B,M\rangle=\sum_{i=1}^n\lambda_i\tilde{m}_{ii}\geq 0$.
\end{proof}
\begin{theorem}
\label{existence}If $A+A^T\prec 0$ then the semidefinite program
defined in Theorem \ref{semidefinite} has an optimal solution.
\end{theorem}
\begin{proof}
First we show that the set $S$ of all matrices $M\succeq 0$
satisfying the equality constraints $\langle
A_i,M\rangle=-p_i(0)$, $i=1,2,\ldots,n-1$, is non-empty. Indeed,
the matrix
\begin{displaymath}
M=\mbox{diag}(-p_1(0)/2a_{11},-p_2(0)/2a_{22},\ldots,-p_n(0)/2a_{nn})
\end{displaymath}
satisfies these conditions and, additionally, it is $M\succeq 0$,
since $p_i(0)\geq 0$ and $a_{ii}<0$ ($A+A^T\prec 0$). Note that
$S$ is closed and convex. Now consider the function
$f:S\rightarrow \Re$ defined by $f(M)=\langle A_n,M\rangle$ and
the matrix $B=-\sum_{i=1}^nA_i=-(A+A^T)\succ 0$. From Lemma
\ref{BM} and for $M\in S$, we have
\begin{displaymath}
\langle B,M\rangle\geq 0\Rightarrow\langle
A_n,M\rangle\leq-\sum_{i=1}^{n-1}\langle
A_i,M\rangle=\sum_{i=1}^{n-1}p_i(0)<+\infty\Rightarrow
f(M)<+\infty\;.
\end{displaymath}
Thus $\sup_{M\in S}f(M)<+\infty$ and since $S$ is closed the
supremum is achieved for a $M_0\in S$, so it is actually a
maximum. The existence of an optimal solution is established.
\end{proof}

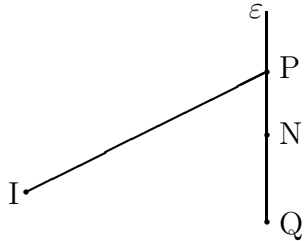
\begin{figure}[h]
\setlength{\unitlength}{0.8cm} \centering
\begin{picture}(3,4) \thicklines
\put(4,0){\line(0,1){3.5}}\put(0,0.5){\line(2,1){4}}
\put(4.2,2.4){P}\put(4,2.5){\circle*{0.1}}\put(4.2,-0.2){Q}\put(4,0){\circle*{0.1}}\put(-0.3,0.3){I}\put(0,0.5){\circle*{0.1}}\put(4.2,1.3){N}\put(4,1.45){\circle*{0.1}}\put(3.7,3.4){$\varepsilon$}
\end{picture}
\caption{\label{reachable}Construction of the reachable set of
point I.}
\end{figure}

We finally show how the semidefinite programming formalism can be
used for calculating the reachable set of point
$I(p_1(0),p_2(0),\ldots,p_n(0))$. Consider the line $\varepsilon$
parallel to $p_n$-axis, with $p_i=\mbox{constant}\geq 0$,
$i=1,2,\ldots,n-1$. The maximum achievable value of $p_n$ on
$\varepsilon$, starting from $I$, can be found by solving the
following semidefinite program: Find $\max_M\langle A_n,M\rangle$
subject to $\langle A_i,M\rangle=p_i-p_i(0)$ for
$i=1,2,\ldots,n-1$ and $M\succeq 0$. If this program has a
solution $M_0$ such that $p_n=p_n(0)+\langle A_n,M_0\rangle\geq
0$, then let P be the point $(p_1,p_2,\ldots,p_n)$ of
$\varepsilon$, see Fig. \ref{reachable}. This point belongs to the
reachable set of $I$. Additionally, every point
$N(p_1,p_2,\ldots,p_{n-1},p_n')$ of $\varepsilon$ with $0\leq
p_n'\leq p_n$, see Fig. \ref{reachable}, belongs also to the
reachable set (first arrive at $P$ and then use
$m=(0,0,\ldots,1)^T$ to go down, since (\ref{psystem}) gives
$\dot{p}_n=a_{nn}<0$, $\dot{p}_i=0$ for $i\neq n$). Thus, the
segment $PQ$, where $Q(p_1,p_2,\ldots,p_{n-1},0)$, is in the
reachable set. By repeating the above procedure for all the
allowed $\varepsilon\parallel p_n$, the reachable set of $I$ can
be constructed.

\section{\label{Ranks}Remarks on the Rank of the Semidefinite Program Solutions}

In the preceding section we saw that the bigger the rank of the
optimal $M$ for the semidefinite program, the more complicated is
the optimal control law. Thus, it would be useful to know if there
exist low rank optimal solutions and, additionally, rank upper
bounds for them. Even more, we would like to know for what
matrices $A$ the corresponding semidefinite program has solutions
of the lowest possible rank $r=1$. In this section we present a
series of results in these directions.
\begin{lemma}
\label{proposition131} Let us fix $A_1,A_2,...,A_k\in\mbox{Sym}_n$
and $\alpha_1,\alpha_2,\ldots,\alpha_k\in \Re$. If there is a
matrix $M\succeq 0$ such that
\begin{displaymath}
\langle A_i,M\rangle=\alpha_i\,,\quad i=1,2,\ldots,k,
\end{displaymath}
then there is a matrix $M_0\succeq 0$ such that
\begin{displaymath}
\langle A_i,M_0\rangle=\alpha_i\,,\quad i=1,2,\ldots,k
\end{displaymath}
and, additionally,
\begin{displaymath}
\mbox{rank}
M_0\leq\left\lfloor\frac{\sqrt{8k+1}-1}{2}\right\rfloor\;,
\end{displaymath}
where $\lfloor\cdot\rfloor$ denotes the integer part of the
embraced number.
\end{lemma}
\begin{proof}
See \cite{Barvinok02}, chapter $II$, proposition $13.1$.
\end{proof}
\begin{proposition}
\label{rankn}If $A+A^T\prec 0$ then there is an optimal solution
$M_0$ to the semidefinite program defined in Theorem
\ref{semidefinite}, with
\begin{displaymath}
\mbox{rank}
M_0\leq\left\lfloor\frac{\sqrt{8n+1}-1}{2}\right\rfloor\;.
\end{displaymath}
\end{proposition}
\begin{proof}
From Theorem \ref{existence} we have that, since $A+A^T\prec 0$,
the semidefinite program has an optimal solution $M\succeq 0$,
which satisfies
\begin{displaymath}
\langle A_i,M\rangle=-p_i(0)\,,\quad
i=1,2,\ldots,n-1\,,\quad\langle A_n,M\rangle=\mathcal{E}\;.
\end{displaymath}
According to Lemma \ref{proposition131}, there exists a
$M_0\succeq 0$ such that
\begin{displaymath}
\langle A_i,M_0\rangle=-p_i(0)\,,\quad
i=1,2,\ldots,n-1\,,\quad\langle A_n,M_0\rangle=\mathcal{E}
\end{displaymath}
and
\begin{displaymath}
\mbox{rank}
M_0\leq\left\lfloor\frac{\sqrt{8n+1}-1}{2}\right\rfloor\;.
\end{displaymath}
Obviously $M_0$ is also an optimal solution.
\end{proof}
\begin{corollary}
\label{rank2} If $A+A^T\prec 0$ and $A$ is $2\times 2$, then the
semidefinite program has an optimal solution of rank $r\leq 1$.
\end{corollary}
\begin{proof}
Apply Proposition \ref{rankn} for $n=2$.
\end{proof}
In general, the bound imposed by Lemma \ref{proposition131} is the
best possible. However, there is one special case where it can be
sharpened.
\begin{lemma}
\label{proposition134} For some positive integer $r$, let us fix
 $k=(r+2)(r+1)/2$ matrices $A_1,A_2,...,A_k\in\mbox{Sym}_n$, where $n\geq r+2$, and $k$ numbers
$\alpha_1,\alpha_2,\ldots,\alpha_k\in \Re$. If there is a matrix
$M\succeq 0$ such that
\begin{displaymath}
\langle A_i,M\rangle=\alpha_i\,,\quad i=1,2,\ldots,k
\end{displaymath}
and the set of all such matrices is bounded, then there is a
matrix $M_0\succeq 0$ such that
\begin{displaymath}
\langle A_i,M_0\rangle=\alpha_i\,,\quad i=1,2,\ldots,k
\end{displaymath}
and, additionally,
\begin{displaymath}
\mbox{rank} M_0\leq r\;.
\end{displaymath}
\end{lemma}
\begin{proof}
See \cite{Barvinok02}, chapter $II$, proposition $13.4$.
\end{proof}
\begin{proposition}
\label{rank3} If $A+A^T\prec 0$ and $A$ is $3\times 3$, then the
semidefinite program has an optimal solution of rank $r\leq 1$.
\end{proposition}
\begin{proof}
Since $A+A^T\prec 0$, the semidefinite program has an optimal
solution $M\succeq 0$, which satisfies
\begin{displaymath}
\langle A_i,M\rangle=-p_i(0)\,,\quad i=1,2\,,\quad\langle
A_n,M\rangle=\mathcal{E}\;.
\end{displaymath}
From Lemma \ref{proposition134} we see that the choice $r=1$ gives
$k=3=n$, since $n=3$ according to the hypothesis of the
proposition. In order to apply Lemma \ref{proposition134}, we just
need to show that the set of optimal matrices, i.e. all the
matrices $M\succeq 0$ satisfying the above relations, is bounded.
Consider the matrix $B=-\sum_{i=1}^3A_i=-(A+A^T)\succ 0$. For a
matrix $M\succeq 0$ in the set of optimal solutions, we have
\begin{displaymath}
\langle B,M\rangle=p_1(0)+p_2(0)-\mathcal{E}<+\infty\;.
\end{displaymath}
But from Lemma \ref{BM}, we have also $\langle
B,M\rangle=\sum_{i=1}^3\lambda_i\tilde{m}_{ii}$, where $\lambda_i$
are the eigenvalues of $B$ and $\tilde{m}_{ii}$ the diagonal
elements of the matrix $\tilde{M}=O^TMO$, $O$ the orthogonal
matrix diagonalizing $B$. Combining these we find that
\begin{displaymath}
\sum_{i=1}^3\lambda_i\tilde{m}_{ii}<+\infty\;.
\end{displaymath}
Since $\lambda_i>0$ and $\tilde{m}_{ii}\geq 0$, the above relation
implies that
\begin{displaymath}
\tilde{m}_{ii}<+\infty\;,
\end{displaymath}
thus
\begin{displaymath}
\mbox{tr}(M)=\mbox{tr}(\tilde{M})<+\infty\;.
\end{displaymath}
But
\begin{displaymath}
\langle M,M\rangle=\mbox{tr}(M^2)\leq (\mbox{tr}(M))^2<+\infty\;,
\end{displaymath}
since $M\succeq 0$. So indeed the set of optimal $M$ is bounded
and we can apply Lemma \ref{proposition134} with $r=1$. This means
that there is an optimal $M_0\succeq 0$ with $\mbox{rank}M_0\leq
1$. Note that the bound that Proposition \ref{rankn} gives in this
case is only $\lfloor(\sqrt{8\cdot 3+1}-1)/2\rfloor=2$.
\end{proof}
\begin{lemma}
\label{problem3} Let us call an $n\times n$ matrix $A=\{a_{ij}\}$
$r$-diagonal if $a_{ij}=0$ unless $|i-j|<r$. Suppose that the
matrices $A_1,A_2,...,A_k\in\mbox{Sym}_n$ are $r$-diagonal and
there exists a matrix $M\succeq 0$ such that
\begin{displaymath}
\langle A_i,M\rangle=\alpha_i\in\Re\,,\quad i=1,2,\ldots,k\;.
\end{displaymath}
Then there exists a matrix $M_0\succeq 0$ such that
\begin{displaymath}
\langle A_i,M_0\rangle=\alpha_i\,,\quad i=1,2,\ldots,k
\end{displaymath}
and, additionally,
\begin{displaymath}
\mbox{rank} M_0\leq r\;.
\end{displaymath}
\end{lemma}
\begin{proof}
See \cite{Barvinok02}, chapter $IV$, corollary $10.3$, problem
$3$.
\end{proof}
\begin{proposition}
If $A+A^T\prec 0$ and $A$ is $r$-diagonal, then the semidefinite
program has an optimal solution of $\mbox{rank}\leq r$.
\end{proposition}
\begin{proof}
Since $A+A^T\prec 0$, there exists an optimal solution of the
semidefinite program. Since $A$ is $r$-diagonal, the corresponding
$A_i$ are also $r$-diagonal. Thus, we can apply Lemma
\ref{problem3}, which assures the existence of an optimal solution
of $\mbox{rank}\leq r$.
\end{proof}

We conclude this section by noting that there is strong numerical
evidence that the following conjecture is true.
\begin{conjecture}
If $A+A^T\prec 0$ then the semidefinite program has an optimal
solution of rank $r=1$.
\end{conjecture}

\section{\label{Examples}Examples}

In this section we solve problems \ref{problem} and \ref{mproblem}
for some specific systems. We start from the system with
\begin{displaymath}
A=\left[\begin{array}{cc}-\xi & -1
\\1 & -\xi\end{array}\right],\;\xi>0\;,
\end{displaymath}
which corresponds to system (\ref{system2}) appeared in the
introduction. It is not necessary to solve numerically the
corresponding semidefinite program, because we can attack this
particular case analytically. Since
$A+A^T=\mbox{diag}(-2\xi,-2\xi)\prec 0$ and $A$ is $2\times 2$,
from Corollary \ref{rank2} we have that there is an optimal
solution to the semidefinite program with rank $r\leq1$. Thus,
there is an optimal constant vector $m=(m_1,m_2)^T$, solution to
Problem \ref{mproblem}. The system equation (\ref{psystem}) with
$A$ given above and $m$ constant gives
\begin{equation}
p_1(T_f)=p_1(0)-(\xi m_1^2+m_1m_2)T_f\,,\quad
p_2(T_f)=p_2(0)+(m_2m_1-\xi m_2^2)T_f\;.
\end{equation}
Optimality requires
\begin{equation}
p_1(T_f)=0\Rightarrow T_f=\frac{p_1(0)}{\xi m_1^2+m_1m_2}\;,
\end{equation}
so
\begin{equation}
p_2(T_f)=p_2(0)+\frac{m_2m_1-\xi m_2^2}{\xi m_1^2+m_1m_2}p_1(0)\;.
\end{equation}
In order to maximize $p_2(T_f)$, we just need to maximize the
coefficient of $p_1(0)$. If we set $m_2/m_1=x$, then this
coefficient takes the form
\begin{equation}
f(x)=\frac{x-\xi x^2}{x+\xi}\;.
\end{equation}
Before maximizing $f$, we find the allowed values of variable $x$.
It should be $p_2(T_f)\geq p_2(0)\Rightarrow x-\xi x^2\geq 0$ and
$p_1(T_f)\leq p_1(0)\Rightarrow x+\xi\geq 0$. These are both
satisfied when $x\in [0,1/\xi]$. We calculate the maximum of $f$
in this interval. It is not difficult to verify that
\begin{equation}
f'(x)=-\frac{\xi(x^2+2\xi x-1)}{(x+\xi)^2}
\end{equation}
becomes zero at the point
\begin{equation}
x_0=\sqrt{1+\xi^2}-\xi
\end{equation}
of the interval $[0,1/\xi]$. Also verify that $f'(x)>0$ for $x\in
[0,x_0)$ and $f'(x)<0$ for $x\in (x_0,1/\xi]$. So $f(x_0)$ is a
maximum in the interval $[0,1/\xi]$. After some manipulation we
find that
\begin{equation}
f(x_0)=x_0^2\;.
\end{equation}
The maximum achievable value of $p_2$ is
\begin{equation}
p_2(T_f)=p_2(0)+x_0^2p_1(0)
\end{equation}
and the optimal unit vector is
\begin{equation}
m=(\frac{1}{\sqrt{1+x_0^2}}\,,\frac{x_0}{\sqrt{1+x_0^2}})\;.
\end{equation}
The optimal trajectory in $p$-space is a straight line joining the
points $(p_1(0),p_2(0))$ and $(0,p_2(T_f))$.

\begin{figure}[p]
\centering
\begin{tabular}{c}
(a)\psfig{file=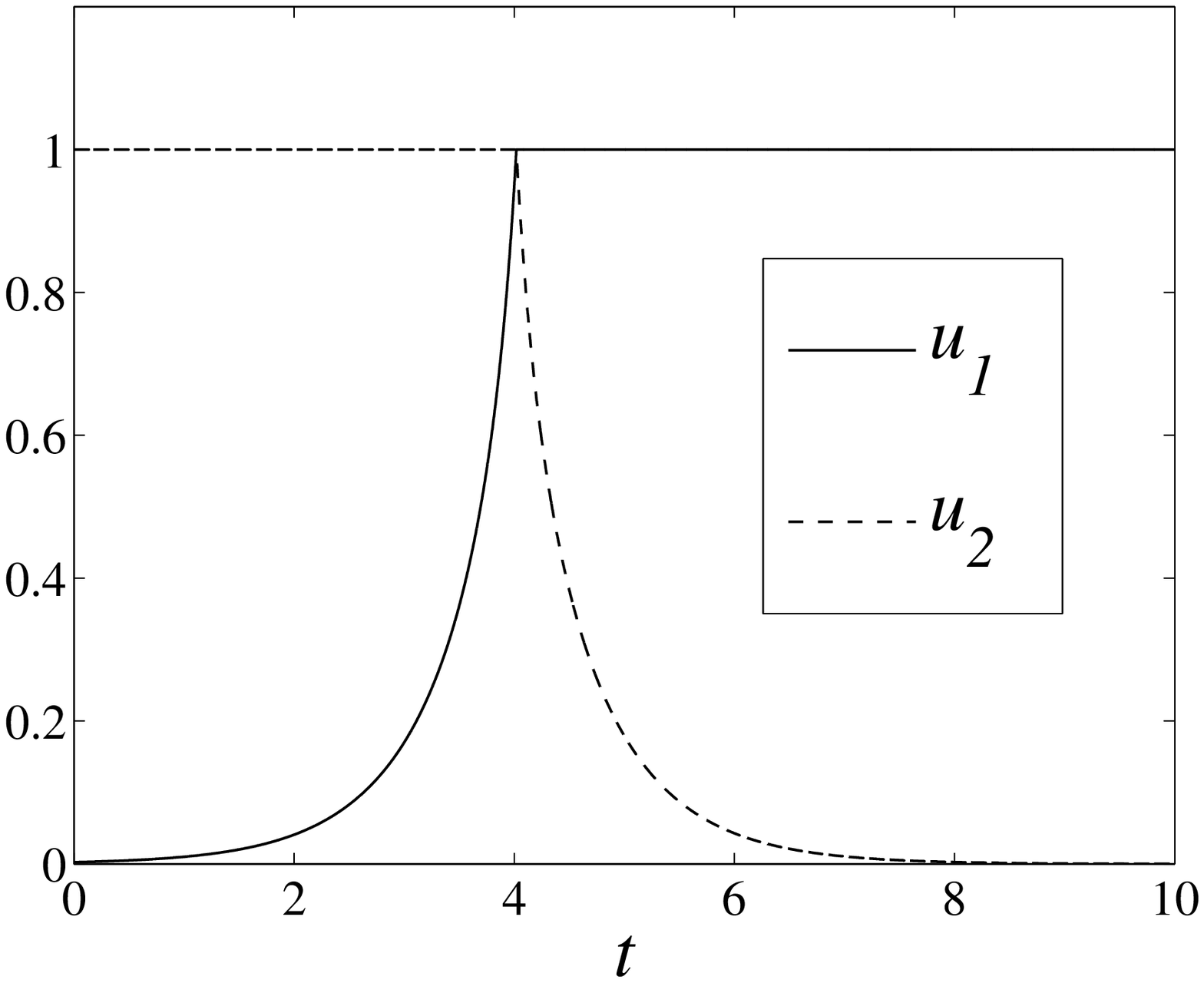,scale=0.35}\\
(b)\psfig{file=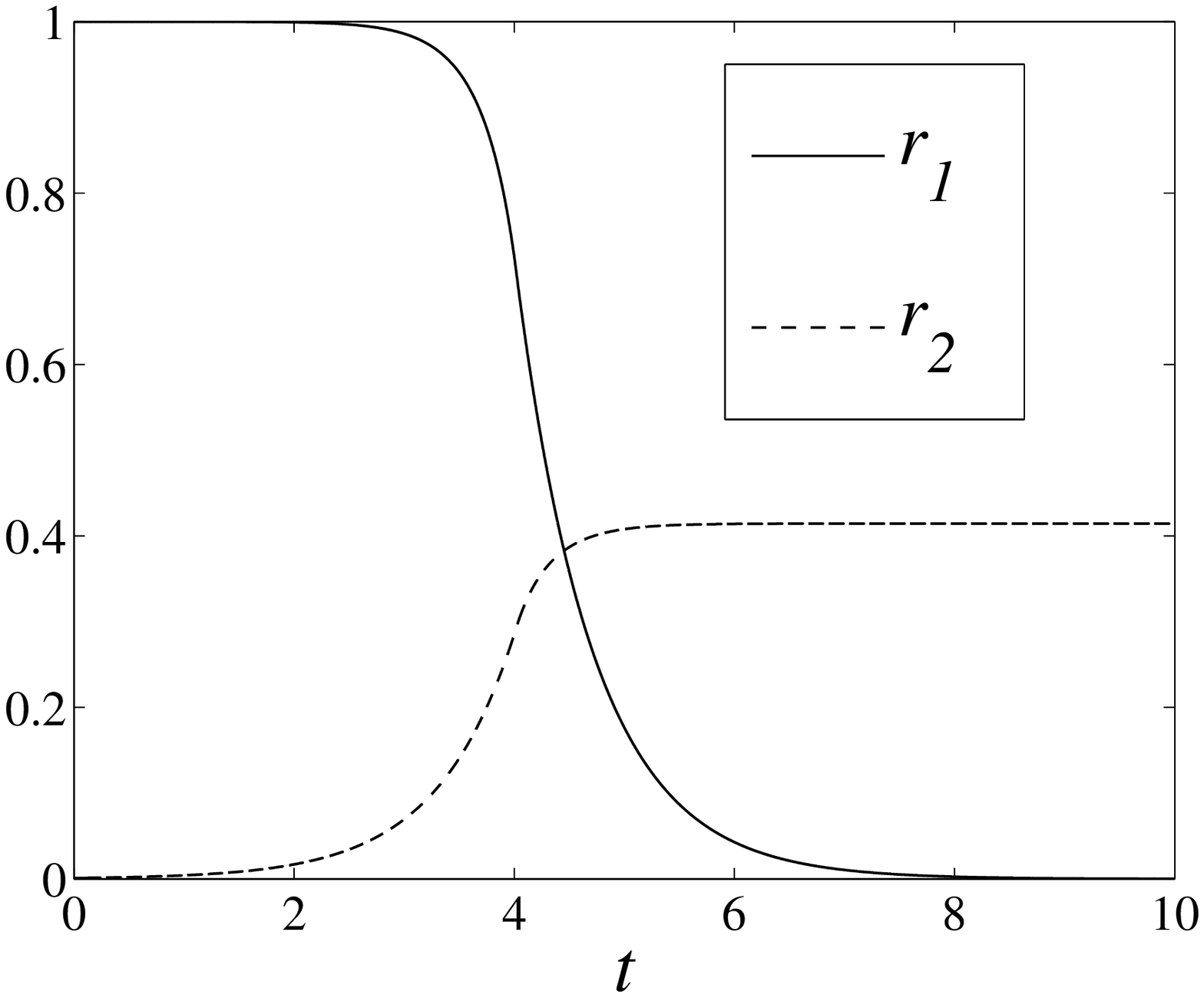,scale=0.35}\\
(c)\psfig{file=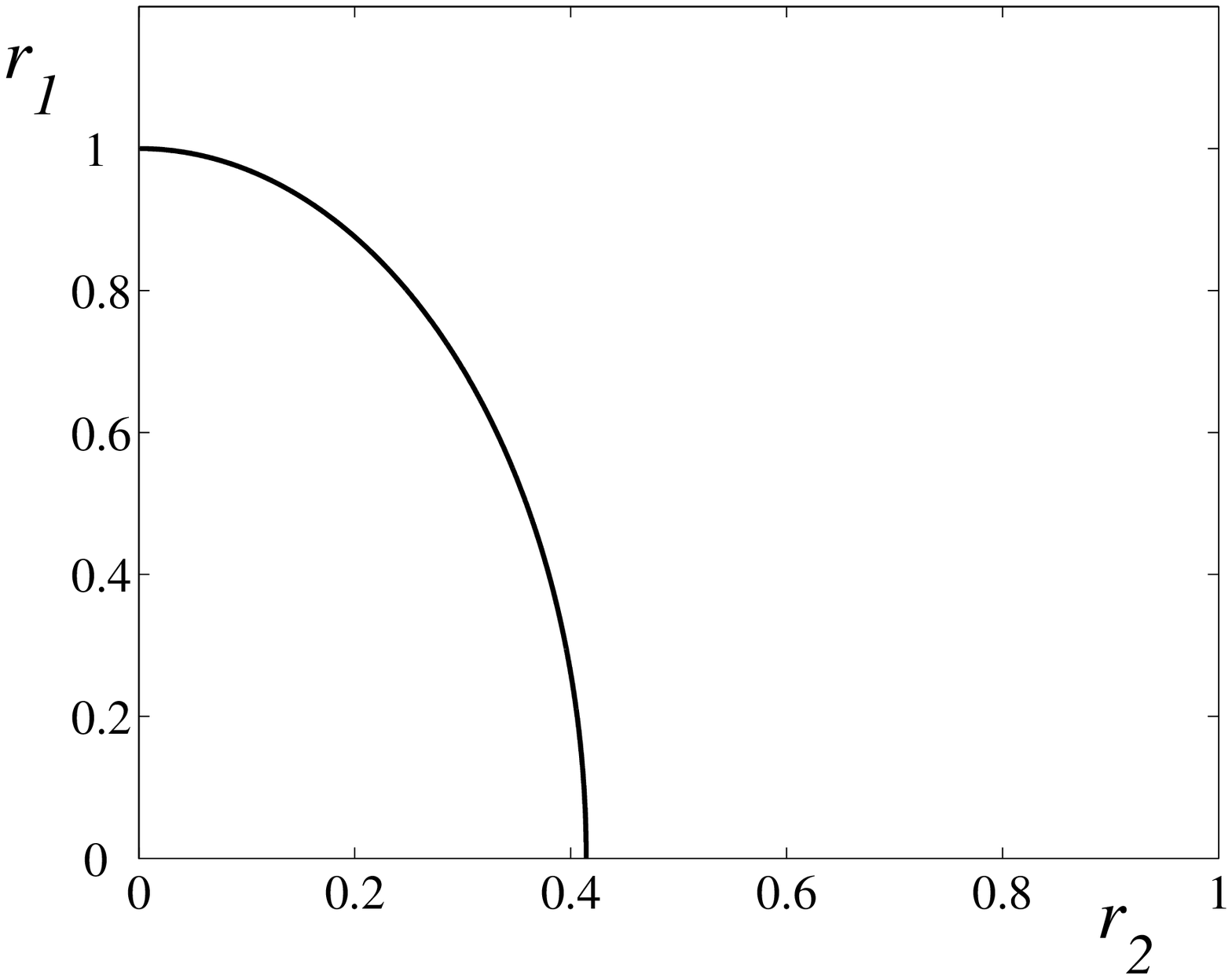,scale=0.35}
\end{tabular}
\caption{\label{2by2}(a) Optimal controls $u_1(t)$ and $u_2(t)$
for system (\ref{system2}) when $\xi=1$ and
$(r_1(0),r_2(0))=(1,0)$ (b) The corresponding state variables
$r_1(t)$ and $r_2(t)$ (c) The optimal trajectory in $r$-space.}
\end{figure}

The maximum achievable value of $r_2$ is
\begin{equation}
r_2(\infty)=\sqrt{r_2^2(0)+x_0^2r_1^2(0)}\;.
\end{equation}
If the starting state is the point $(r_1(0), r_2(0))=(1,0)$, the
maximum transfer efficiency takes the value
\begin{equation}
r_2(\infty)=x_0=\sqrt{1+\xi^2}-\xi\;.
\end{equation}
For $\xi=1$ this efficiency is $\sqrt{2}-1$. The optimal controls
$u_1, u_2$ for system (\ref{system2}), can be found by using the
method described in section \ref{SolStrat}. If we define
\begin{displaymath}
\mathcal{M}=\max\left( 1, \frac{x_0r_1}{r_2} \right)\;,
\end{displaymath}
the optimal policy can be realized as
\begin{displaymath}
u_1=\frac{1}{\mathcal{M}},\quad\quad u_2 =\frac{x_0r_1}{r_2}u_1\;.
\end{displaymath}
Observe that the initial point $(1,0)$ is a stationary point of
the optimal control policy [$r_2(0)=0\Rightarrow
\mathcal{M}=\infty\Rightarrow u_1=0\Rightarrow u_2=0$]. This
optimal policy in the infinite horizon case should then be
interpreted as the limit of optimal control policy for the
corresponding finite time problem [Finite time for the $r_i$
system, don't confuse it with the finite time problem for $p_i$
which corresponds to the infinite horizon problem for $r_i$. It is
the special case $k_1=k_2=0$ of the finite time problem solved in
the preceding chapter. The solution for this particular case can
be found in \cite{Khaneja03}]. In practice, we give a small but
finite value in $r_2(0)$ (an initial 'kick' from zero) which makes
the optimal control law applicable. In Fig. \ref{2by2}(a) we plot
the optimal controls $u_1$ and $u_2$. In Fig. \ref{2by2}(b) we
depict $r_1(t), r_2(t)$ and in Fig. \ref{2by2}(c) the
corresponding optimal trajectory in $r$-space. For all these
figures it is $\xi=1$ and $(r_1(0), r_2(0))=(1,0)$.

\begin{remark}
The closure of the reachable set of point $(1,0)$ is
\begin{displaymath}
\overline{\textbf{R}((1,0))}=\{r_1,r_2\geq
0\mid\sqrt{r_2^2+x_0^2r_1^2}\leq x_0\}\;,
\end{displaymath}
where $x_0=\sqrt{1+\xi^2}-\xi$. This set is depicted in Fig.
\ref{2by2}(c) for $\xi=1$. The closure of the reachable set
$\overline{\textbf{R}((1,0,0,0))}$ for the corresponding bilinear
system (\ref{physicalsystem}) is
\begin{displaymath}
\{(x_1,x_2,y_1,y_2)\in
\Re^4\mid\sqrt{(x_2^2+y_2^2)+x_0^2(x_1^2+y_1^2)}\leq x_0\}\;.
\end{displaymath}
\end{remark}

The next case that we examine is the system with
\begin{displaymath}
A=\left[\begin{array}{ccc}-\xi & -1 & 0 \\
1 & -\xi & -1\\
0 & 1 & -\xi \end{array}\right],\;\xi>0\;.
\end{displaymath}
Since $A+A^T=\mbox{diag}(-2\xi,-2\xi,-2\xi)\prec 0$ and $A$ is
$3\times 3$, from Proposition \ref{rank3} we have that the
semidefinite program has a solution of rank $r\leq 1$. Now let us
become more specific, so set $\xi=1$ and consider the starting
point $(p_1(0),p_2(0),p_3(0))=(1,1,0)$. The corresponding matrices
$A_i$ are
\begin{displaymath}
A_1=\left[\begin{array}{ccc}-2 & -1 & 0 \\-1 & 0 & 0\\
0 & 0 & 0 \end{array}\right]\,,\quad A_2=\left[\begin{array}{ccc}0 & 1 & 0 \\1 & -2 & -1\\
0 & -1 & 0 \end{array}\right]\,,\quad A_3=\left[\begin{array}{ccc}0 & 0 & 0 \\0 & 0 & 1\\
0 & 1 & -2 \end{array}\right]\;.
\end{displaymath}
If we solve numerically the corresponding semidefinite program
using some appropriate software package, for example SDPT3
\cite{Toh99}, we find that the optimal matrix $M\succeq 0$ is
\begin{displaymath}
M=\left[\begin{array}{ccc}0.1775 & 0.3225 & 0.1304 \\0.3225 & 0.5856 & 0.2368\\
0.1304 & 0.2368 & 0.0958 \end{array}\right]\;,
\end{displaymath}
and the maximum achievable value of $p_3$ is
\begin{displaymath}
p_3(T_f)=p_3(0)+\langle A_3,M\rangle=\langle
A_3,M\rangle=0.2821\;.
\end{displaymath}
It is easy to verify that this matrix has two zero eigenvalues and
one nonzero, so its rank is indeed $r=1$. It can be written in the
form $M=\lambda mm^T$, where $\lambda=0.8589\,(=T_f)$ is the
nonzero eigenvalue and
\begin{displaymath}
m=(m_1,m_2,m_3)^T=(0.4546,0.8257,0.3339)^T
\end{displaymath}
the corresponding eigenvector. This unit vector is the optimal
solution for Problem \ref{mproblem}.

The maximum achievable value of $r_3$ is
$r_3(\infty)=\sqrt{p_3(T_f)}=0.5311$. We find the optimal
$u_1,u_2,u_3$. Let us set
\begin{displaymath}
x_0=\frac{m_2}{m_1}=1.8163\,,\quad y_0=\frac{m_3}{m_1}=0.7345\;.
\end{displaymath}
If we define
\begin{displaymath}
\mathcal{M}=\max(1,\frac{x_0r_1}{r_2}, \frac{y_0r_1}{r_3})\;,
\end{displaymath}
then the optimal policy can be realized as
\begin{displaymath}
u_1=\frac{1}{\mathcal{M}},\quad u_2=\frac{x_0r_1}{r_2}u_1,\quad
u_3=\frac{y_0r_1}{r_3}u_1\;.
\end{displaymath}
\begin{figure}[p]
\centering
\begin{tabular}{c}
(a)\psfig{file=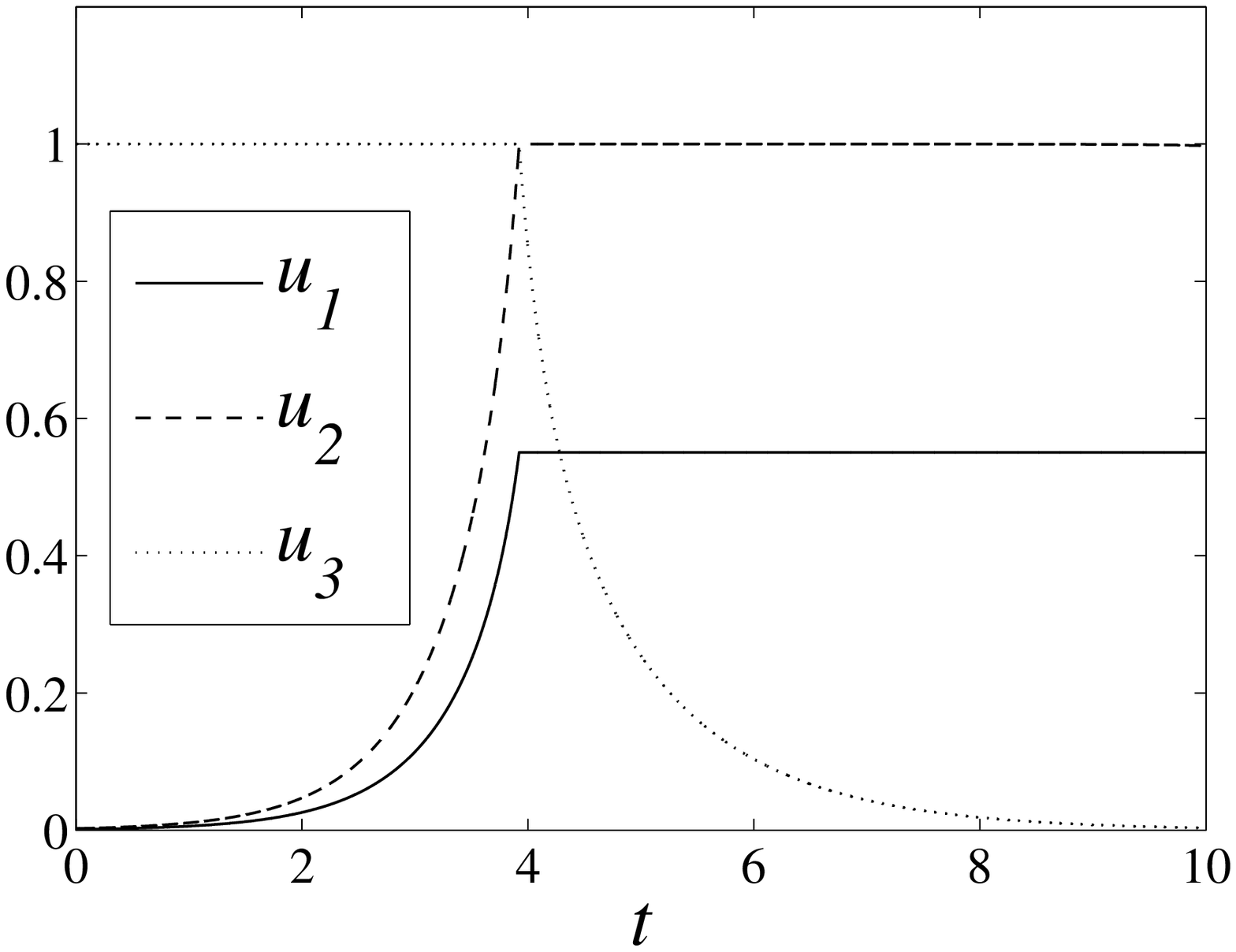,scale=0.35}\\
(b)\psfig{file=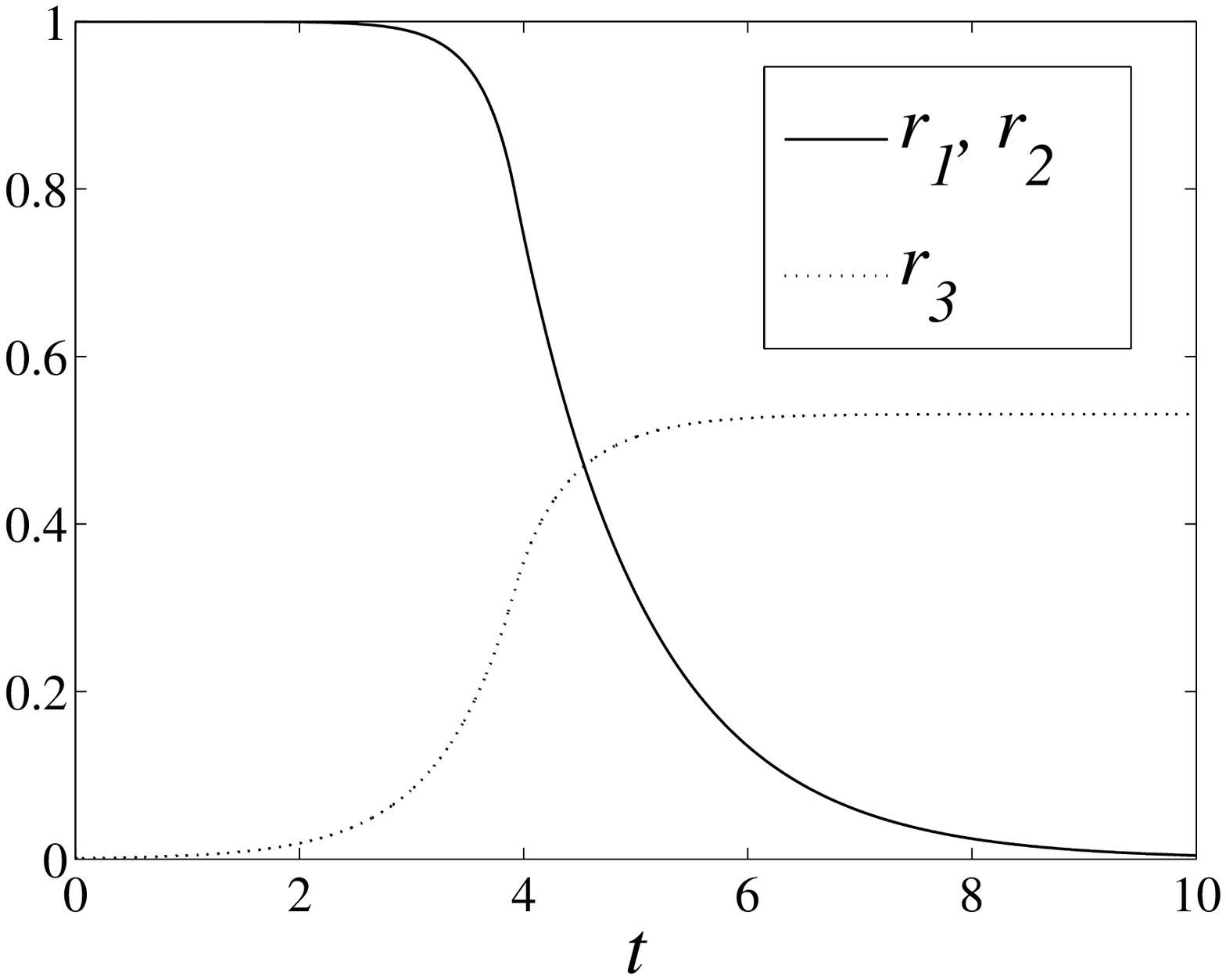,scale=0.35}\\
(c)\psfig{file=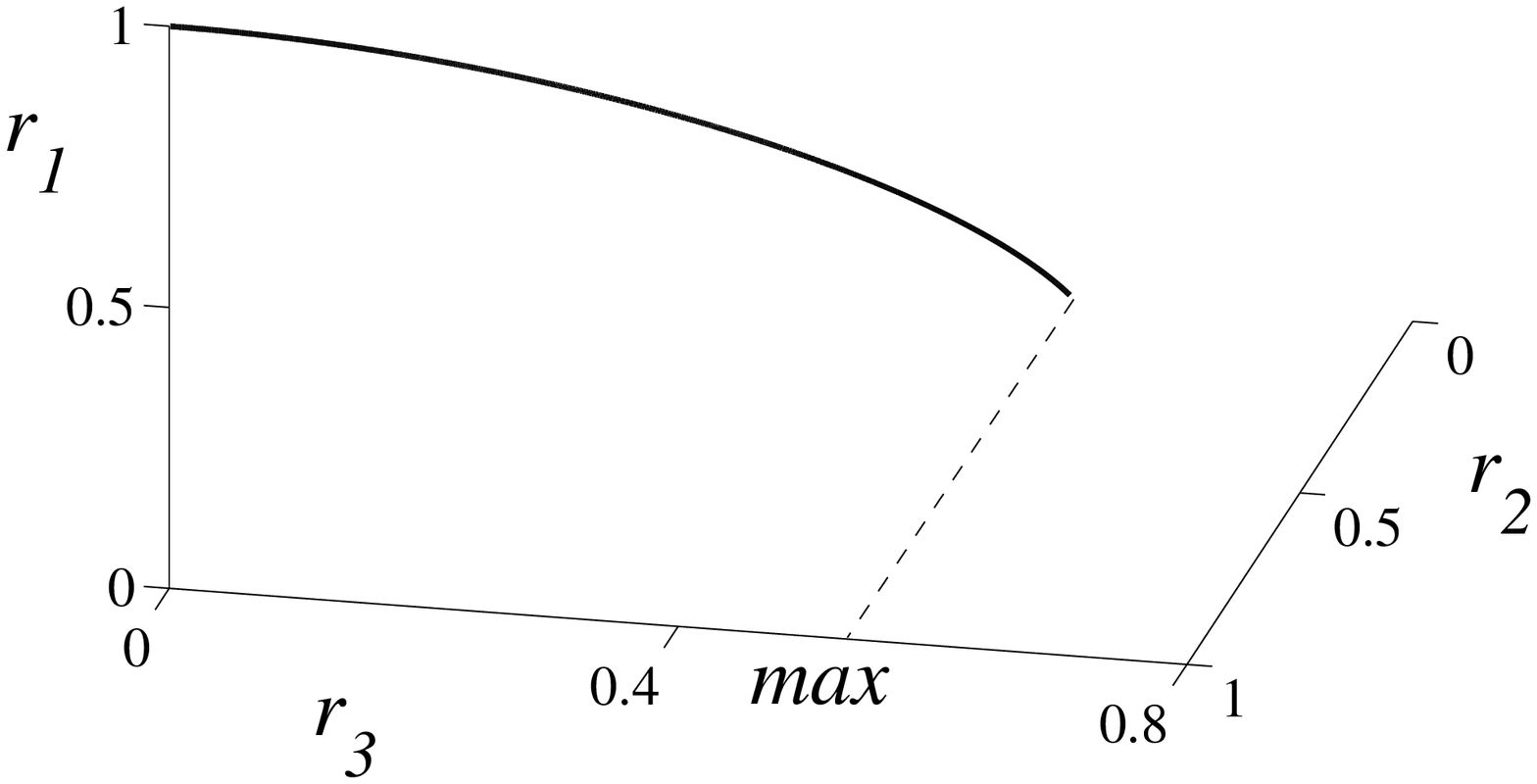,scale=0.35}
\end{tabular}
\caption{\label{r3}(a) Optimal controls $u_1(t),u_2(t)$ and
$u_3(t)$ for system (\ref{generalsystem}), with $A$ the $3\times
3$ matrix given in the text, when $\xi=1$ and
$(r_1(0),r_2(0),r_3(0))=(1,1,0)$ (b) The corresponding state
variables $r_1(t),r_2(t)$ and $r_3(t)$. Observe that
$r_2(t)/r_1(t)=1$ throughout. Remember that the optimal trajectory
in $p$-space is a straight line ending at the point
$(0,0,p_3(T_f))$, so $p_2(t)/p_1(t)=p_2(0)/p_1(0)=1$ for the
starting point $(1,1,0)$ (c) The optimal trajectory in $r$-space.}
\end{figure}
Observe that the initial point $(1,1,0)$ is a stationary point of
the optimal control policy [$r_3(0)=0\Rightarrow
\mathcal{M}=\infty\Rightarrow u_1=0\Rightarrow u_2=u_3=0$]. The
situation is similar with that in the previous example. Again, the
optimal policy in the infinite horizon case should be interpreted
as the limit of the optimal policy for the corresponding finite
time problem. In practice, a small but finite value is given to
$r_3(0)$. In Fig. \ref{r3}(a) we plot the optimal controls
$u_1,u_2$ and $u_3$. In Fig. \ref{r3}(b) we depict $r_1(t),
r_2(t), r_3(t)$ and in Fig. \ref{r3}(c) the corresponding optimal
trajectory in $r$-space.

Another interesting case to examine is the same system with
$\xi>0$ unspecified and starting point $(p_1(0),0,p_3(0))$. This
problem can be solved analytically and has the practical
application that it gives an upper bound for our ability to
coherently control a specific dissipative quantum system
\cite{Stefanatos}. As before, we know that there is an optimal
constant vector $m=(m_1,m_2,m_3)^T$. From equations
(\ref{psystem}) we find
\begin{eqnarray}
p_1(T_f) & = & p_1(0)-(\xi m_1^2+m_1m_2)T_f\;,\nonumber\\
p_2(T_f) & = & p_2(0)+(m_2m_1-\xi
m_2^2-m_2m_3)T_f\;,\\
p_3(T_f) & = & p_3(0)+(m_3m_2-\xi m_3^2)T_f\nonumber\;.
\end{eqnarray}
Optimality requires
\begin{equation}
p_1(T_f)=0\Rightarrow T_f=\frac{p_1(0)}{\xi m_1^2+m_1m_2}
\end{equation}
and
\begin{equation}
\label{p2zero} p_2(T_f)=0\Rightarrow m_2m_1-\xi m_2^2-m_2m_3=0\;.
\end{equation}
So, we have to maximize
\begin{equation}
p_3(T_f)=p_3(0)+\frac{m_3m_2-\xi m_3^2}{\xi m_1^2+m_1m_2}p_1(0)\;
\end{equation}
subject to the constraint (\ref{p2zero}). We just need to maximize
the coefficient of $p_1(0)$ under the same condition. If we set
\begin{displaymath}
\frac{m_2}{m_1}=x\,,\quad\frac{m_3}{m_1}=y\;,
\end{displaymath}
then this coefficient takes the form
\begin{equation}
g(x,y)=\frac{xy-\xi y^2}{x+\xi}\;,
\end{equation}
while the condition becomes
\begin{equation}
\label{condition} x(1-\xi x-y)=0\Rightarrow y=1-\xi x\;.
\end{equation} Note that
x=0 gives $g\leq 0$ so it is rejected. Using (\ref{condition}),
$g$ becomes a function of $x$ only
\begin{equation}
f(x)=g(x,y(x))=\frac{-\xi(1+\xi^2)x^2+(1+2\xi^2)x-\xi}{x+\xi}\;.
\end{equation}
We find the allowed values of $x$. A natural requirement is
$p_3(T_f)\geq p_3(0)\Rightarrow y(x-\xi y)\geq 0\Rightarrow (y\geq
0\;\mbox{and}\;x-\xi y\geq 0)\;\mbox{or}\;(y\leq
0\;\mbox{and}\;x-\xi y\leq 0)$. Using (\ref{condition}) we find
that the first option implies $x_1\leq x\leq x_2$ and the second
$x_2\leq x\leq x_1$, where $x_1=\xi/(1+\xi^2),x_2=1/\xi$. Since
$x_1<x_2$, only the first option is acceptable, so it must be
$x\in [x_1,x_2]$. For such $x$, the similar requirement
$p_1(T_f)\leq p_1(0)\Rightarrow x+\xi\geq 0$ is satisfied. So we
maximize $f$ in the interval $[x_1,x_2]=[\xi/(1+\xi^2),1/\xi]$.
Solving the equation $f'(x_0)=0$, we find
\begin{equation}
x_0=\sqrt{\xi^2+2}-\xi.
\end{equation}
Indeed, $x_0\in [x_1,x_2]$. The corresponding maximum value of $f$
is
\begin{equation}
f_{max}=\frac{x_0^4}{4}\;.
\end{equation}
The maximum  achievable value of $p_3$ is
\begin{equation}
p_3(T_f)=p_3(0)+f_{max}p_1(0)\;.
\end{equation}

\begin{figure}[p]
\centering
\begin{tabular}{c}
(a)\psfig{file=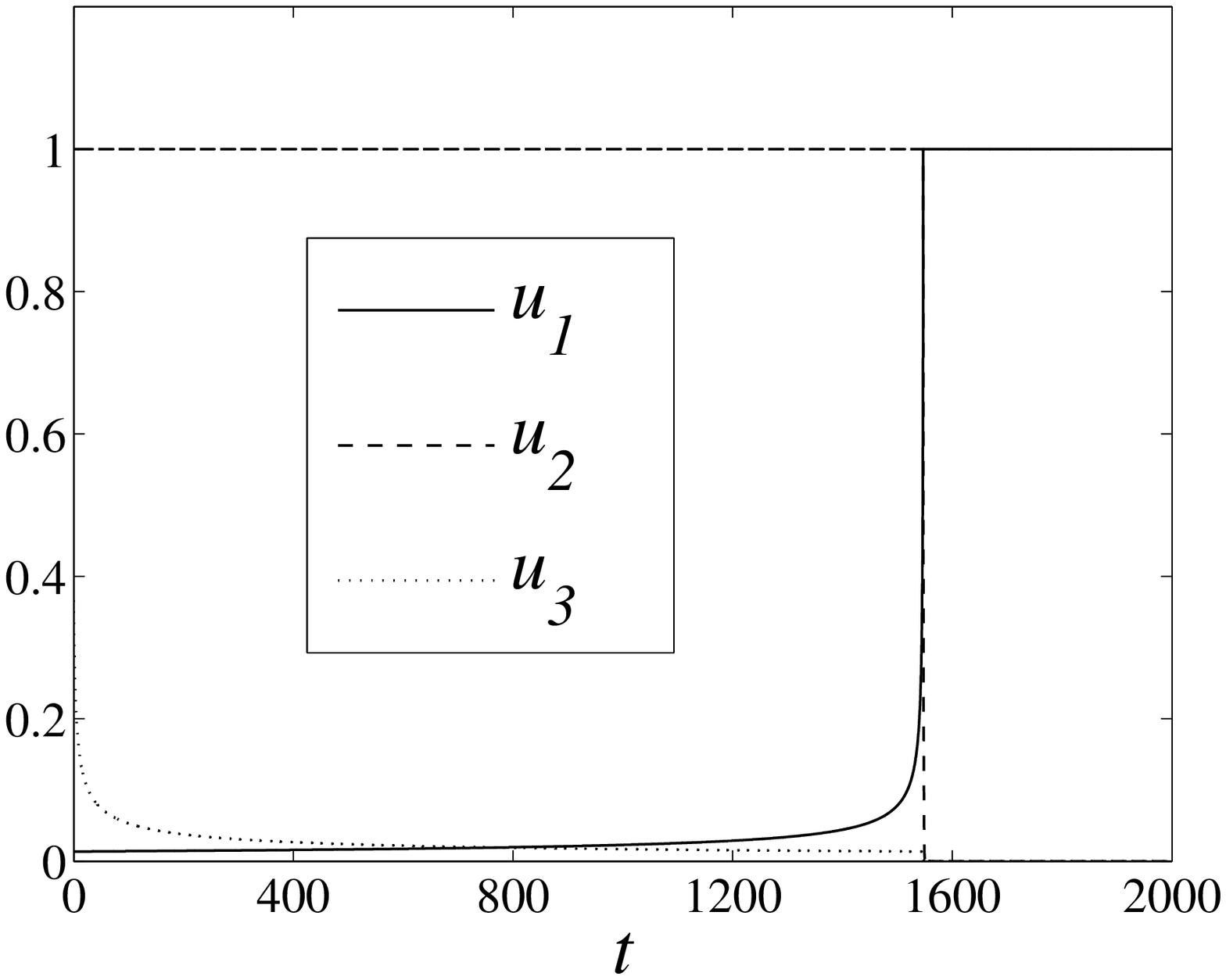,scale=0.35}\\
(b)\psfig{file=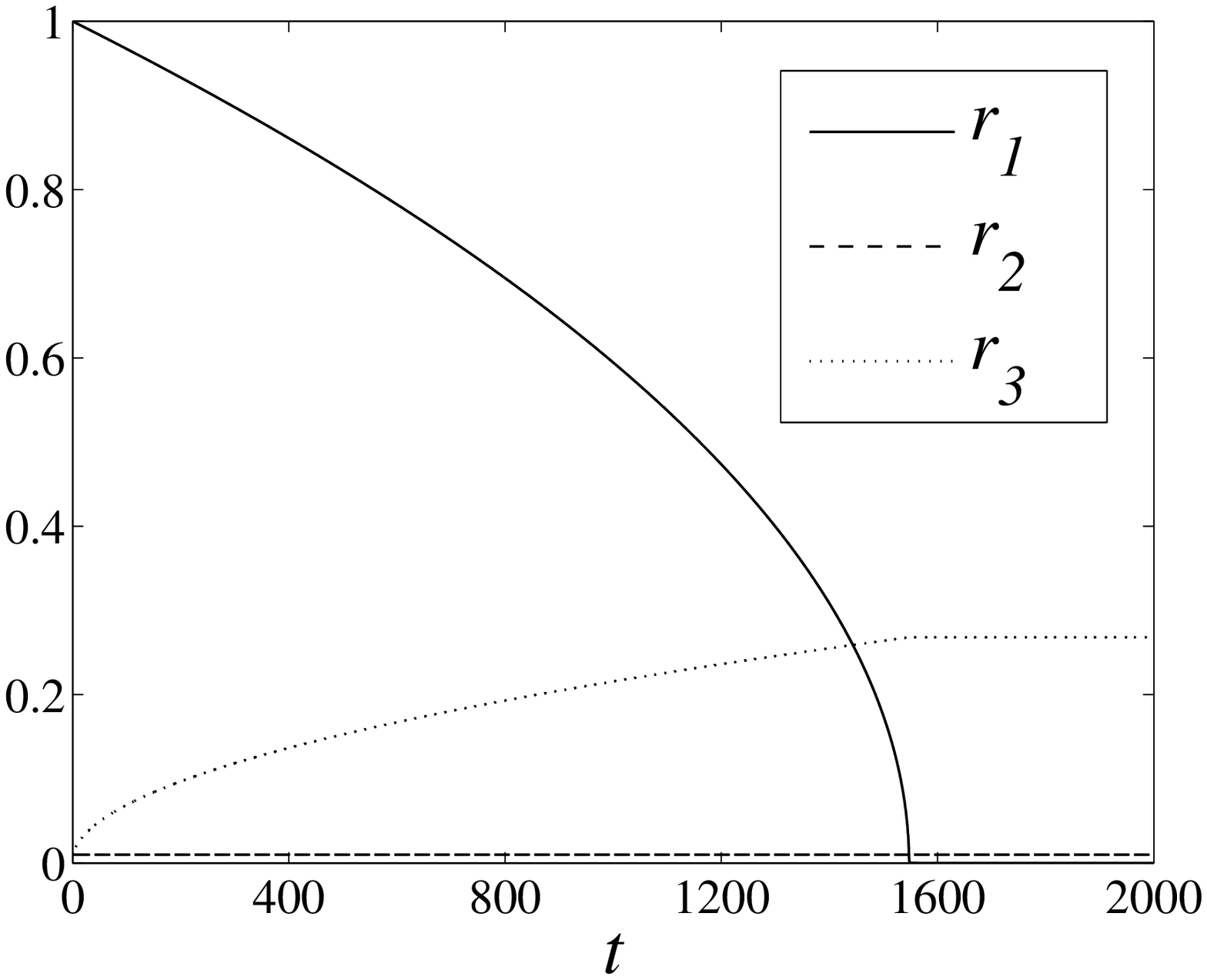,scale=0.35}\\
(c)\psfig{file=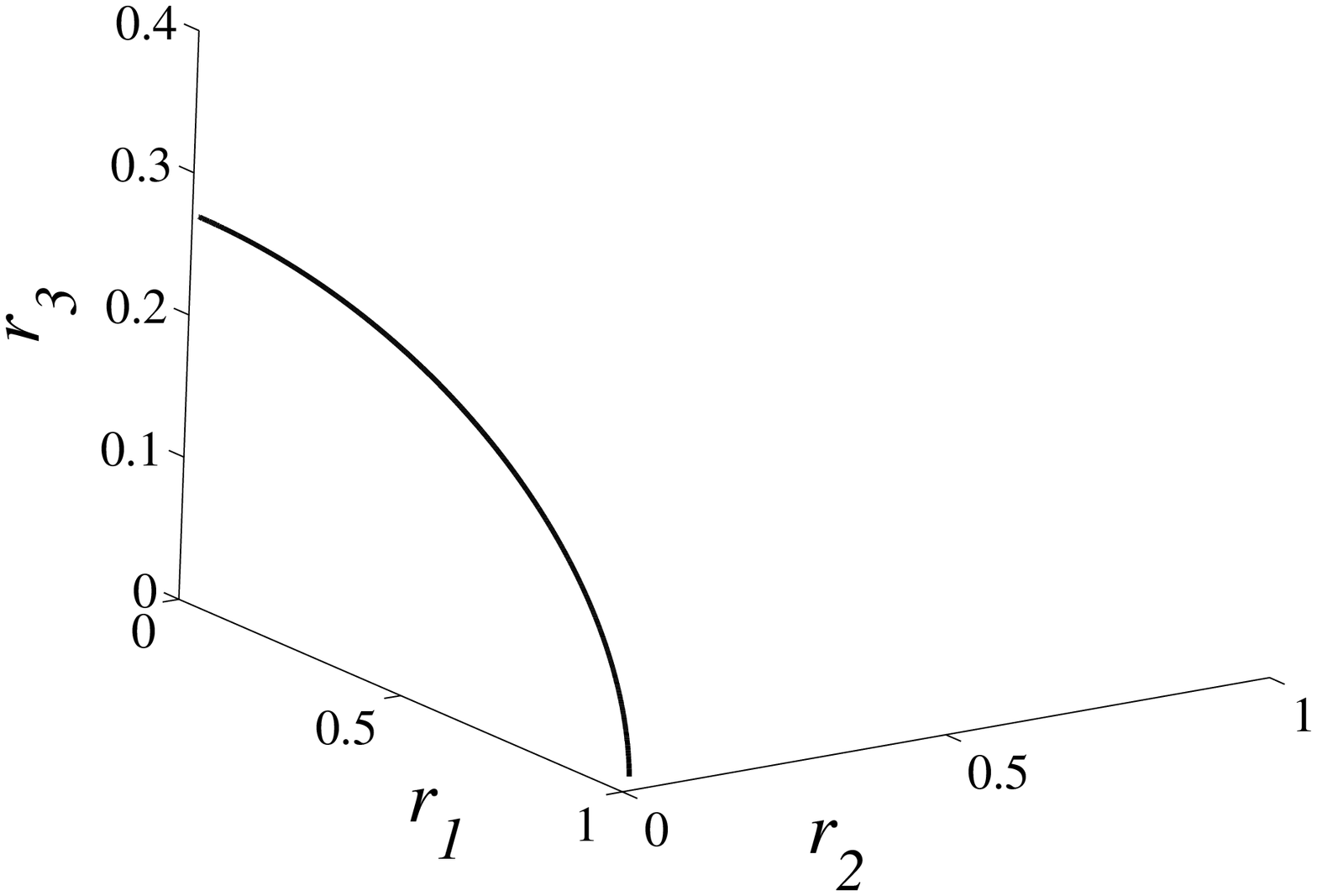,scale=0.35}
\end{tabular}
\caption{\label{rinfty}(a) Optimal controls $u_1(t),u_2(t)$ and
$u_3(t)$ for system (\ref{generalsystem}), with $A$ the $3\times
3$ matrix given in the text, when $\xi=1$ and
$(r_1(0),r_2(0),r_3(0))=(1,\epsilon,\epsilon)$, $0<\epsilon\ll 1$.
Here we take $\epsilon=0.01$ for convenience (b) The corresponding
state variables $r_1(t),r_2(t)$ and $r_3(t)$. Note that transfer
$r_1\rightarrow r_3$ takes place through $r_2$ which is held to
the small constant value $r_2=\epsilon$. Thus, this transfer
requires more time compared to the preceding examples (c) The
optimal trajectory in $r$-space.}
\end{figure}

Condition (\ref{p2zero}) implies that in the optimal case it is
$\dot{p}_2=0$, so it is also $\dot{r}_2=0$. If $r_2(0)=0$ then
$r_2(t)=0$ and, as we can see from (\ref{generalsystem}), there is
no transfer from $r_1$ to $r_3$. What we actually examine here is
the limiting case $r_2(0)=\epsilon\rightarrow0^+$, where
$\epsilon$ is an arbitrarily small positive number. We can still
use condition (\ref{p2zero}), i.e. $\dot{r}_2=0$. The transfer
$r_1\rightarrow r_3$ takes place through $r_2$ which is held to
the small constant value $r_2=\epsilon$. The maximum achievable
value of $r_3$, which corresponds to the limit
$\epsilon\rightarrow0^+$, is
\begin{equation}
r_3(\infty)=\sqrt{r_3^2(0)+f_{max}r_1^2(0)}.
\end{equation}
If the starting state is the point $(1,\epsilon,0)$, where
$\epsilon\rightarrow0^+$, the maximum efficiency is
\begin{equation}
\label{effi}
r_3(\infty)=\sqrt{f_{max}}=\frac{x_0^2}{2}=\frac{(\sqrt{\xi^2+2}-\xi)^2}{2}.
\end{equation}
For $\xi=1$ we find that this efficiency is $2-\sqrt{3}$. In Fig.
\ref{rinfty} we plot the optimal controls $u_1(t),u_2(t),u_3(t)$,
the state variables $r_1(t),r_2(t),r_3(t)$ and the optimal
trajectory in r-space. Observe that the starting point is actually
$(1,\epsilon,\epsilon)$. It is necessary to give a small positive
initial value to $r_3$, since the point $(1,\epsilon,0)$ is still
a stationary point of the optimal policy. If the starting point is
$(1,\epsilon,\epsilon)$, then by solving the corresponding
semidefinite program we find numerically the same efficiency as in
(\ref{effi}), in the limit $\epsilon\rightarrow0^+$.

\section{\label{Con} Conclusion}

In this paper we studied a class of bilinear control systems,
motivated by optimal control problems arising in the context of
dissipative quantum dynamics. It was shown that the optimal
solution and the reachable set of these systems can be found by
solving a semidefinite program. As a practical result, solutions
to these problems give upper bounds for the ability to coherently
control quantum mechanical phenomena in presence of dissipation.
In the area of coherent spectroscopy, these results translate into
the maximum signal that can be obtained in an experiment. The
paper also motivates the use of semidefinite programming to study
reachable sets of more general bilinear control systems.

\end{document}